\def \CC {\mathbb C}

\def \RR {\mathbb R}
\def \ZZ {\mathbb Z}

\def \epsilon{\varepsilon}

\def \d {\text{d}}

\def \bfmu {{\boldsymbol{\mu}}}
\def \bflambda {{\boldsymbol{\lambda}}}

\def \fine {{\hfill \qedsymbol}}

\def \si {\sigma}

\renewcommand{\S}{{\mathcal S}}

\documentclass[12pt,reqno]{amsart}

\usepackage{amsmath,amssymb}

\begin{document}


\title[]{Some remarks on the convergence of the Dirichlet series of $L$-functions and related questions}
\author[]{J.KACZOROWSKI \lowercase{and} A.PERELLI}
\maketitle

{\bf Abstract.} First we show that the abscissae of uniform and absolute convergence of Dirichlet series coincide in the case of $L$-functions from the Selberg class $\S$. We also study the latter abscissa inside the extended Selberg class, indicating a different behavior in the two classes. Next we address two questions about majorants of functions in $\S$, showing links with the distribution of the zeros and with independence results.



\vskip.5cm
\section{Introduction}

\smallskip
Let 
\[
F(s)=\sum_{n=1}^\infty \frac{a(n)}{n^s}
\]
be a Dirichlet series which converges somewhere in the complex plane. It is well known that there are four classical abscissae associated with $F(s)$: the abscissa of {\it convergence} $\si_c(F)$, of {\it uniform convergence} $\si_u(F)$, of {\it absolute convergence} $\si_a(F)$ and of {\it boundedness} $\si_b(F)$. It may well be that $\si_c(F)=-\infty$, in which case the other three abscissae equal $-\infty$ as well. From the theory of Dirichlet series we know that
\[
\sigma_c(F)\leq \sigma_b(F)= \sigma_u(F)\leq \sigma_a(F),
\]
and in general this is best possible, i.e. inequalities cannot be replaced by equalities. We refer to Maurizi-Queff\'elec \cite{Ma-Qu/2010} for a modern reference for this sort of problems.

\medskip
Our first result is that $\si_b(F)=\si_a(F)$ for an important class of Dirichlet series, namely those defining the $L$-functions of the {\it Selberg class} $\S$. We recall that the axiomatic class $\S$ contains, at least conjecturally, most $L$-functions from number theory and automorphic forms theory, and that $\si_b(F)=\si_a(F)$ is known in some special cases like the Riemann or the Dedekind zeta functions. The Selberg class $\S$ is defined, roughly, as the class of Dirichlet series absolutely convergent for $\si>1$, having analytic continuation to $\CC$ with at most a pole at $s=1$, satisfying a functional equation of Riemann type and having an Euler product representation. Moreover, their coefficients satisfy the Ramanujan condition $a(n)\ll n^\epsilon$ for any $\epsilon>0$. We also recall that the {\it extended Selberg class} $\S^\sharp$ is the larger class obtained by dropping the Euler product and Ramanujan condition requirements in the definition of $\S$. We refer to our survey papers \cite{Ka-Pe/1999b}, \cite{Kac/2006}, \cite{Per/2005}, \cite{Per/2004}, \cite{Per/2010} and to the forthcoming book \cite{Ka-Pe/book} for definitions, examples and the basic theory of the Selberg classes $\S$ and $\S^\sharp$. In particular, we refer to these papers for the definition of {\it degree} $d_F$, {\it conductor} $q_F$ and {\it standard twist} of $F(s)$.

\medskip
{\bf Theorem 1.} {\sl Suppose that $F(s)$ belongs to the Selberg class. Then}
\begin{equation}
\label{1.1}
\sigma_b(F)= \sigma_u(F) = \sigma_a(F).
\end{equation}

\medskip
Several months after submitting this result, the note by Brevig-Heap \cite{Br-He/2015} appeared, where the authors prove the same theorem in the much more general framework of Dirichlet series with multiplicative coefficients. Trying to understand Brevig-Heap's proof, based on Bohr's theory, we noticed that their result was already known to Bohr himself in 1913 (see \cite{Boh/1913}, Satz XI, p.480); incidentally, Bohr's paper  \cite{Boh/1913} appears as item [5] of the reference list in Brevig-Heap \cite{Br-He/2015}. We wish to thank Dr. Mattia Righetti for bringing \cite{Br-He/2015} and \cite{Boh/1913} to our attention and for his advice concerning these papers. We decided to keep Theorem 1 since our proof is  different, easier and more direct; moreover, some points in the proof will be useful for the other results in the paper.

\medskip
We expect that actually $\sigma_a(F)=1$ {\it for all} $F\in\S$. This is known for most classical $L$-functions and, in the general case of the class $\S$, under the assumption of the {\it Selberg orthonormality conjecture}; however, an unconditional proof is missing at present. See again the above quoted references for definitions and results about such a conjecture. 

\medskip
Note that the abscissa of convergence $\si_c(F)$ can be smaller than 1 for functions in $\S$. For example, the Dirichlet $L$-functions $L(s,\chi)$ with a primitive non-principal character $\chi$ are convergent in the half-plane $\si>0$. Actually, several general results are known about the abscissa $\si_c(F)$ for functions $F(s)$ in the extended Selberg class $\S^\sharp$. First of all
\begin{equation}
\label{1.2}
\text{{\it if $F\in\S^\sharp$ is entire with degree $d\geq1$, then \ $\frac12 - \frac{1}{2d} \leq \si_c(F) \leq 1-\frac{2}{d+1}$}}
\end{equation}
(recall that {\it there exist no functions $F\in\S^\sharp$ with degree} $0<d<1$, see \cite{Ka-Pe/1999a} and Conrey-Ghosh \cite{Co-Gh/1993}). Indeed, the first inequality in \eqref{1.2} is Corollary 3 in \cite{Ka-Pe/2005} and is based on the properties of the standard twist, while the second inequality follows from a well known theorem of Landau \cite{Lan/1915}. 
Moreover, in accordance with classical degree 2 conjectures and with the general $\Omega$-theorem in Corollary 2 of \cite{Ka-Pe/2005}, we expect that equality holds in the left inequality in \eqref{1.2}. Further
\[
\text{{\it $\si_c(F)=-\infty$ if and only if $d_F=0$,}}
\]
since the degree 0 functions of $\S^\sharp$ are Dirichlet polynomials (see \cite{Ka-Pe/1999a}). From \eqref{1.2} we also deduce that
\[
\text{{\it $\si_c(F)=1$ if and only if $F(s)$ has a pole at $s=1$.}}
\]
We also remark that {\it if the Lindel\"of Hypothesis holds for $F\in\S^\sharp$, then $\si_c(F)\leq 1/2$}.

\medskip
The behavior of $\si_a(F)$ in the extended class $\S^\sharp$ is different from the expected behavior in $\S$. Indeed, in the next section, which is also of independent interest, we show that 
\[
\text{{\it there exist functions $F\in\S^\sharp$ with $\si_a(F)$ arbitrarily close to} 1/2.}
\] 
We conclude this section with the following

\medskip
{\bf Question.} Does \eqref{1.1} hold for the functions in the extended Selberg class ? \fine

\medskip
A variant of the question is: does \eqref{1.1} hold for linear combinations
\[
F(s) = \sum_{j=1}^N c_jF_j(s)
\]
with $F_j\in\S$ and $c_j\in\CC$ ? If needed, one may assume that $F(s)$ belongs to $\S^\sharp$.

\medskip
Since $\sigma_a(F)=1$ for most classical $L$-functions $F(s)$, Theorem 1 prevents the possibility of getting information on the non-trivial zeros exploiting the properties of the abscissa of uniform convergence. On the other hand, if $F\in\S$ is bounded for $\si>1-\delta$ for some $\delta>0$, then its Dirichlet series is absolutely convergent for $\si>1-\delta$ and hence $F(s)\neq0$ by Euler's identity. In the next theorems we replace boundedness by more general majorants and deduce some consequences.

\medskip
Let $F\in\S$ be of degree $d$, $N_F(\si,T)$ be the number of zeros $\rho=\beta+i\gamma$ with $\beta>\si$ and $|\gamma|\leq T$, and denote the {\it density abscissa} $\si_D(F)$ by
\[
\si_D(F) = \inf\{\si: N_F(\si,T)=o(T)\}.
\]
An inspection of the proof of Lemma 3 in \cite{Ka-Pe/2003}, obtained by a rudimentary version of Montgomery's zero-detecting method, shows that
\[
N_F(\si,T) \ll T^{4(d+3)(1-\si)+\epsilon}.
\]
Hence in general
\[
\frac12 \leq \si_D(F) \leq 1-\frac{1}{4(d+3)},
\]
although it is well known that the classical $L$-functions $F(s)$ of degree 1 and 2 have $\si_D(F)=1/2$, see e.g. Luo \cite{Luo/1995}. Actually, one can prove that $\si_D(F)=1/2$ for all $F\in\S$ with degree $0<d\leq 2$. Further, let $f(s)$ be holomorphic in $\si>1-\delta$ for some $\delta>0$ and almost periodic on the line $\si=A$ for some $A>1$. We say that $f(s)$ is a $\delta$-{\it almost periodic majorant} of $F(s)$ if
\begin{equation}
\label{1.3}
|F(s)|\leq c(\si)|f(s)|
\end{equation}
in the half-plane $\si>1-\delta$, where $c(\si)>0$ is a continuous function for $\si>1-\delta$.

\medskip
{\bf Theorem 2.} {\sl Let $F\in\S$ and $f(s)$ be a $\delta$-almost periodic majorant of $F(s)$. Then $F(s)$ and $f(s)$ have the same zeros, with the same multiplicity, in the half-plane $\si>\max(1-\delta,\si_D(F))$.}

\medskip
{\bf Remark.} Clearly, in view of \eqref{1.3} each zero of $f(s)$ is also a zero of $F(s)$; the non-trivial part of Theorem 2 says that the opposite assertion holds true as well. Note that we do not require that $f(s)$ is almost periodic for $\si>1-\delta$, but only on some vertical line far on the right. We already noticed that, as a consequence of Theorem 1, $F(s)\neq0$ in every right half-plane where it is bounded. An immediate consequence of Theorem 2 is that $F(s)\neq0$ for $\si>\max(1-\delta,\si_D(F))$ if $f(s)$ is a non-vanishing  $\delta$-almost periodic majorant.  In particular, from the density estimates reported above when $d\leq 2$, if $\delta=1/2$ then the Riemann Hypothesis holds for such $F(s)$. \fine

\medskip
Our final result is a kind of new independence statement for $L$-functions from the Selberg class. Several forms of independence are known in $\S$, such as the linear independence, the multiplicity one property and the orthogonality conjecture and some of its consequences; see our above quoted surveys on the Selberg class. The new independence result is expressed in terms of majorants as follows.

\medskip
{\bf Theorem 3.} {\sl Let $F,G\in\S$ be such that $F(s)\ll |G(s)|$ for $\si>1/2$. Then $F(s)=G(s)$.}

\medskip
The special nature of the majorant is very important here. Indeed, suppose that $G(s)$ is entire; then Theorem 2 gives only that $F(s)$ and $G(s)$ have the same zeros for $\si>\si_D(F)$. Instead, exploiting the information that $G\in\S$, Theorem 3 shows that actually $F(s)=G(s)$. In other words, no function from $\S$ can dominate in $\si>1/2$ another function from $\S$. We may regard this as a weak form of a well known result obtained, under stronger assumptions, by Selberg \cite{Sel/1992} and Bombieri-Hejhal \cite{Bo-He/1995} about the statistical independence of the values of $L$-functions.

\medskip
{\bf Acknowledgements.} 
This research was partially supported by Istituto Nazionale di Alta Matematica and by grant 2013/11/B/ST1/02799 {\sl Analytic Methods in Arithmetic} of the National Science Centre.

\section{The lift operator}

\smallskip
Let $Q>0$, $\bflambda=(\lambda_1,\dots,\lambda_r)$ with $\lambda_j>0$, $\bfmu=(\mu_1,\dots,\mu_r)$ with $\mu_j\in\CC$ and $|\omega|=1$. We denote by $W(Q,\bflambda,\bfmu,\omega)$  the $\RR$-linear space of the Dirichlet series solutions $F(s)$ of the functional equation
\begin{equation}
\label{1.4}
Q^s \prod_{j=1}^r\Gamma(\lambda_js+\mu_j)F(s) = \omega Q^{1-s} \prod_{j=1}^r\Gamma(\lambda_j(1-s)+\overline{\mu_j})\overline{F(1-\overline{s})}.
\end{equation}
Given an integer $k\geq 1$, we define the {\it $k$-lift operator} by
\[
F(s)\longmapsto F_k(s)=F(ks+\frac{1-k}{2});
\]
clearly, the operator is trivial for $k=1$. A simple computation shows that 
\begin{equation}
\label{1.5}
\text{{\it if $F\in W(Q,\bflambda,\bfmu,\omega)$ then $F_k\in W(Q^k,k\bflambda,\bfmu + \frac{1-k}{2}\bf\bflambda,\omega)$.}}
\end{equation}
In particular, from \eqref{1.5} we have that degree $d_{F_k}$ and conductor $q_{F_k}$ of $F_k(s)$ satisfy
\begin{equation}
\label{1.6}
d_{F_k} = kd_F \hskip1.5cm  q_{F_k} = q_F^k k^{kd_F}.
\end{equation}

We recall (see the above references) that the class $\S^\sharp$ consists of the Dirichlet series satisfying a functional equation of type \eqref{1.4}, where now $\Re{\mu_j}\geq0$, with the following properties: $F(s)$ is absolutely convergent for $\sigma>1$ and $(s-1)^mF(s)$ is entire of finite order for some integer $m\geq0$. Therefore we consider
\[
B_F= 2\min_{1\leq j\leq r} \frac{\Re{\mu_j}}{\lambda_j}+1,
\]
which {\it is an invariant of }$\S^\sharp$ (see again the above references) since a simple computation shows that
\[
B_F=-2\max_{\rho} \Re{\rho} + 1,
\]
where $\rho$ runs over the trivial zeros of $F(s)$. From the definition of the $k$-lift operator and \eqref{1.5} we see that, given $F\in\S^\sharp$, {\it the lifted function $F_k(s)$ also belongs to $\S^\sharp$ provided $1\leq k\leq B_F$ and, if $B_F\geq 2$, $F(s)$ is entire}. Indeed, if $k\geq 2$, $F(s)$ has to be holomorphic at $s=1$ otherwise the pole of $F_k(s)$ is not at $s=1$, and the bound $k\leq B_F$ is needed to have non-negative real part of the $\mu$'s data of $F_k(s)$. Therefore, defining $V(Q,\bflambda,\bfmu,\omega)$ to be the $\RR$-linear space of the entire functions $F\in\S^\sharp$ satisfying \eqref{1.4} (again with $\Re{\mu_j}\geq0$), we have that
\[
\text{{\it for $1\leq k\leq B_F$, the $k$-lift operator maps $V(Q,\bflambda,\bfmu,\omega)$ into $V(Q^k,k\bflambda,\bfmu + \frac{1-k}{2}\bf\bflambda,\omega)$}}.
\]
Note that $B_F$ depends only on $\bflambda$ and $\bfmu$, so it is the same for all functions in $V(Q,\bflambda,\bfmu,\omega)$. Note also that {\it the Selberg class $\S$ is not preserved under the above mappings} since the Ramanujan condition is not (necessarily) satisfied by $F_k(s)$ even if $F(s)$ does; see the examples below. Further, a simple computation shows that {\it the $k$-lift operator commutes with the map sending $F(s)$ to its standard twist}. We also remark that the requirement $\Re{\mu_j}\geq0$ in the definition of $\S^\sharp$, which is responsible for the limitation $k\leq B_F$ in \eqref{1.6}, is apparently not of primary importance in the theory of the Selberg class. Hence, although formally not belonging to $\S^\sharp$, the lifts $F_k(s)$ of entire $F\in\S^\sharp$ with $k>B_F$ are further examples of Dirichlet series with continuation over $\CC$ and functional equation. A similar remark applies to the other condition in the definition of $V(Q,\bflambda,\bfmu,\omega)$, namely the holomorphy at $s=1$.

\medskip
{\bf Examples.} The Riemann zeta function $\zeta(s)$ cannot be lifted inside $\S^\sharp$ since it has $B_\zeta=1$. The same holds for the Dirichlet $L$-functions with even primitive characters, while those associated with odd primitive characters may be lifted inside $\S^\sharp$ for $k=2$ and $k=3$. However, after lifting their Dirichlet coefficients do not satisfy the Ramanujan condition, hence the lifted Dirichlet $L$-functions do not belong to $\S$. Note that, once suitably normalized, the lifts with $k=2$ become the $L$-functions associated with half-integral weight modular forms; see the books by Hecke \cite{Hec/1983} and Ogg \cite{Ogg/1969}. Concerning degree 2, we consider the $L$-functions associated with holomorphic eigenforms of level $N$ and integral weight $K$; see Ogg \cite{Ogg/1969}. Denoting by $F(s)$ their normalization satisfying a functional equation reflecting $s\mapsto1-s$ (instead of the original $s\mapsto K-s$), we have that 
\[
B_F=K.
\]
In other words, the normalized $L$-functions of eigenforms of weight $K$ may be lifted inside $\S^\sharp$ with $k$ up to their weight. Here we consider only eigenforms since in general the $L$-functions of modular forms of level $N$ satisfy a slightly different functional equation, not of $\S^\sharp$ type. \fine

\medskip
We finally turn to the problem of the absolute convergence abscissa in $\S^\sharp$. Let $F\in\S^\sharp$ be of degree $d\geq 1$. Then, thanks again to the properties of the standard twist, we know that 
\begin{equation}
\label{1.7}
\si_a(F) \geq \frac12 + \frac{1}{2d};
\end{equation}
this folows from Theorem 1 of \cite{Ka-Pe/2005}. On the other hand, if $F\in\S^\sharp$ we have that the series
\[
\sum_{n=1}^\infty \frac{|a(n)|}{n^{k\si +(1-k)/2}}
\]
converges for $\si>1/2 + 1/(2k)$. Hence from \eqref{1.6} and \eqref{1.7} we obtain that if both $F(s)$ of degree $d\geq 1$ and $F_k(s)$ belong to $\S^\sharp$, then
\begin{equation}
\label{1.8}
\frac12 + \frac{1}{2kd} \leq \si_a(F_k) \leq \frac12 +\frac{1}{2k}.
\end{equation}
Since the above examples show that there exist functions $F\in\S^\sharp$ with arbitrarily large $B_F$ (e.g. the holomorphic eigenforms with arbitrarily large weight $K$), \eqref{1.8} shows that $\si_a(F)$ can be arbitrarily close to $1/2$ inside $\S^\sharp$. Hence the behavior of $\si_a(F)$ in the extended class $\S^\sharp$ is definitely different from its expected behavior in the class $\S$.

\section{Proof of Theorem 1}

\smallskip
Observe that the case $d=0$ is trivial, since $F(s)$ is identically $1$; see Conrey-Ghosh \cite{Co-Gh/1993}. For $d$ positive we have $\sigma_b(F)\geq 1/2$, since $F(s)$ is unbounded for $\sigma<1/2$ by the functional equation and the properties of the $\Gamma$ function.  Therefore, to prove the assertion it suffices to show the following fact: if for a certain $1/2<\sigma_0\leq 1$ the function $F(s)$ is bounded for $\sigma>\sigma_0$, then $\sigma_a(F)\leq \sigma_0$. 

\medskip
Let us fix an $\varepsilon\in(0,\sigma_0-1/2)$, and let $c_0=c_0(\varepsilon)$ be such that $|a(n)|\leq c_0n^{\varepsilon\slash 2}$ for all $n\geq 1$. Without loss of generality we may assume that $c_0\geq 3$. Consider the finite set of primes
\[
S_{\varepsilon}=\{p: |a(p)|>p^{\varepsilon\slash 2} \  \ \text{or}\ \ p< c_0^{2\slash\varepsilon}\}.
\]
Let 
\begin{equation}
\label{1.9}
F_p(s) = \sum_{m=0}^\infty \frac{a(p^m)}{p^{ms}}
\end{equation}
denote the $p$-th Euler factors of $F(s)$. We split the Euler product as
\begin{equation}
\label{1.10}
\begin{split}
F(s)&=\prod_{p\not\in S_{\varepsilon}}\left(1+\frac{a_F(p)}{p^s}\right) \prod_{p\in S_{\varepsilon}}F_p(s)
\prod_{p\not\in S_{\varepsilon}}\left(F_p(s)\left(1+\frac{a_F(p)}{p^s}\right)^{-1}\right)\\
&=P_1(s)P_2(s)P_3(s),
\end{split}
\end{equation}
say. Both $P_2(s)$ and its inverse $1\slash P_2(s)$ have Dirichlet series representations which  converge absolutely for $\sigma>\theta$ for some $\theta<1\slash 2$. This is a simple consequence of the definition of the Selberg class; see the above quoted references. Therefore, $P_2(s)$ and $1\slash P_2(s)$ are bounded for $\sigma>\si_0$.

\medskip
In view of \eqref{1.9} we have
\[
P_3(s)=\prod_{p\not\in S_{\varepsilon}}\left(1+\sum_{m=2}^{\infty}\frac{b(p^m)}{p^{ms}}\right)
\]
with
\[
b(p^m)=\sum_{l=0}^m(-1)^la(p)^la(p^{m-l}).
\]
Hence, recalling that $p\not\in S_{\varepsilon}$, $m\geq 2$  and $c_0\geq 3$, we have
\[
|b(p^m)|\leq \sum_{l=0}^m|a(p)|^l|a(p^{m-l})|\leq c_0mp^{m\varepsilon\slash 2}\leq p^{m\varepsilon}.
\]
Thus for $\sigma>1/2+\varepsilon$ and $p\not\in S_{\varepsilon}$ we have
\[ 
\sum_{m=2}^{\infty}\frac{|b(p^m)|}{p^{m\sigma}} <1 \quad \text{and} \quad \sum_{p\not\in S_{\varepsilon}} \sum_{m=2}^{\infty}\frac{|b(p^m)|}{p^{m\sigma}}\ll 1.
\]
Hence both $P_3(s)$ and  $1\slash P_3(s)$ are bounded and have Dirichlet series representations which  converge absolutely for $\sigma>\sigma_0$ (recall that $\sigma_0>1/2+\varepsilon$).
 
 \medskip
We therefore see that $P_1(s)=F(s)\slash(P_2(s)P_3(s))$ is bounded for $\sigma>\sigma_0$. Let us write
\[
P_1(s)=\sum_{n=1}^{\infty}\frac{c(n)}{n^s}.
\]
The coefficients $c(n)$ are completely multiplicative, and the series converges for $\sigma>\sigma_0$. Fix such a $\sigma$, and a positive  $\delta<\sigma-\sigma_0$. Consider the following familiar Mellin's transform
 \[ 
 \sum_{n=1}^{\infty}\frac{c(n)}{n^{\sigma+it}}e^{-n\slash Y} = \frac{1}{2\pi i}\int_{1-i\infty}^{1+i\infty} \frac{F(w+\sigma+it)}{P_2(w+\sigma+it)P_3(w+\sigma+it)} \Gamma(w) Y^w \d w.
 \]
 We shift the line of integration to $\Re(w)=-\delta$ and obtain
 \[ 
 \sum_{n=1}^{\infty}\frac{c(n)}{n^{\sigma+it}}e^{-n\slash Y} = \frac{F(\sigma+it)}{P_2(\sigma+it)P_3(\sigma+it)} + O(Y^{-\delta})\ll 1
 \]
 uniformly in $t\in {\mathbb R}$ and $Y\geq 1$. Since $|c(n)|\leq n^{\varepsilon/2}$, due to the decay of the exponential we may cut the sum on the left hand side to $n\leq 3Y\log Y$, say, producing an extra error term of size $O(1\slash Y)$. Thus
 \begin{equation}
 \label{1.11}
  \sum_{n\leq 3Y\log Y}\frac{c(n)}{n^{\sigma+it}}e^{-n\slash Y} \ll 1
 \end{equation} 
 uniformly in $t\in\RR$ and $Y\geq1$.
 
\medskip
Now we apply Kronecker's theorem in the following form, see Theorem 8 of Ch.VIII of Chandrasekharan \cite{Cha/1968}. If $\theta_1,\dots,\theta_k\in\RR$ are linearly independent over $\ZZ$, $\beta_1,\dots,\beta_k\in\RR$ and $T,\eta>0$, then there exist $t>T$ and $n_1,\dots,n_k\in\ZZ$ such that
\begin{equation}
\label{1.12}
|t\theta_\ell - n_\ell -\beta_\ell|<\eta \hskip1.5cm \ell=1,\dots,k.
\end{equation}
We choose the $\theta$'s as $-\frac{1}{2\pi}\log p$ with the primes $p\leq 3Y\log Y$ not in $S_{\varepsilon}$ and, correspondingly, the $\beta$'s such that $|c(p)|=c(p)e^{2\pi i\beta_p}$ for each such $p$. Hence by \eqref{1.12} there exists a sequence of real numbers $t_\nu\to+\infty$ such that
\[
c(p)p^{-it_\nu} \to |c(p)| \hskip1.5cm \nu\to\infty
\]
uniformly for the primes $p\leq 3Y\log Y$ not in $S_\varepsilon$. By the complete multiplicativity of $c(n)$ we infer that
\[ 
c(n)n^{-it_{\nu}} \to |c(n)|  \hskip1.5cm \nu\to\infty
\]
uniformly for $n\leq 3Y\log Y$. Thus putting $t=t_{\nu}$ in (\ref{1.11}) and making $\nu\to\infty$ we obtain
\[ 
\sum_{n\leq Y}\frac{|c(n)|}{n^{\sigma}}\leq e\sum_{n\leq 3Y\log Y}\frac{|c(n)|}{n^{\sigma}}e^{-n\slash Y}
 = e\lim_{\nu\to\infty} \sum_{n\leq 3Y\log Y}\frac{c(n)}{n^{\sigma+it_{\nu}}}e^{-n\slash Y} \ll 1
\]
 uniformly for $Y\geq 1$. Letting $Y\to\infty$, we see that the Dirichlet series of $P_1(s)$ converges absolutely for $\sigma>\sigma_0$. 
 
\medskip
Summarizing, we have shown that the Dirichlet series of $P_1(s)$, $P_2(s)$ and $P_3(s)$ are absolutely convergent for $\si>\si_0$, hence the Dirichlet series of $F(s)$ is also absolutely convergent for $\si>\si_0$ thanks to \eqref{1.10}, and the result follows.

\section{Proof of Theorem 2}

\smallskip
As in Theorem 1 the case $d=0$ is trivial, hence we assume $d>0$ and consider the function
\[
h(s) = \frac{F(s)}{f(s)}
\]
for $\si>1-\delta$. From \eqref{1.3} we have that $h(s)$ is holomorphic for $\si>1-\delta$, bounded on every closed vertical strip inside $\si>1-\delta$ and almost periodic on the line $\si=A$. For a given $\epsilon>0$, let $\tau$ be an $\epsilon$-almost period of $h(A+it)$, namely for every $t\in\RR$
\[
|h(A+i(t+\tau)) - h(A+it)| < \epsilon.
\]
Then, by the convexity following from Phragm\'en-Lindel\"of's theorem applied to $h(s+i\tau)-h(s)$, given $\eta>1-\delta$ and any $\eta<\si<A$ we have
\[
\begin{split}
\sup_{t\in\RR} |h(\si+i(t+\tau)) - h(\si+it)&| \leq \big(\sup_{t\in\RR} |h(\eta+i(t+\tau)) - h(\eta+it)|\big)^{\frac{A-\si}{A-\eta}} \\
&\times \big(\sup_{t\in\RR} |h(A+i(t+\tau)) - h(A+it)|\big)^{\frac{\si-\eta}{A-\eta}}.
\end{split}
\]
Hence we obtain that
\[
\sup_{t\in\RR} |h(\si+i(t+\tau)) - h(\si+it)| \ll \epsilon^c
\]
uniformly in any closed strip contained in $\eta<\si<A$, where $c>0$ depends on the strip. Since $\epsilon$ is arbitrarily small, $h(s)$ is uniformly almost periodic in such strips. Suppose now that $h(\rho)=0$ for some $\rho$ with $\Re{\rho}>1-\delta$. Then by a well known argument based on Rouch\'e's theorem we have that for any $1-\delta<\eta<\Re{\rho}$
\[
T \ll N_h(\eta,T) \leq N_F(\eta,T) =o(T)
\]
if $\eta>\si_D(F)$, a contradiction. Thus $h(s)\neq0$ for $\si>\max(1-\delta,\si_D(F))$, hence every zero of $F(s)$ in this half-plane is a zero of $f(s)$. Theorem 2 is therefore proved, since the opposite implication is a trivial consequence of \eqref{1.3}.

\section{Proof of Theorem 3}

\smallskip
Again the case $d=0$ is trivial, since in this case $F(s)\equiv1$ and so $G(s)$ does not vanish inside the critical strip, thus its degree is 0 and hence $G(s)\equiv1$ as well. Let $F,G\in\S$ be with positive degrees and coefficients $a_F(n)$ and $a_G(n)$, respectively, and consider the function
\[
H(s)=\frac{F(s)}{G(s)} = \sum_{n=1}^\infty \frac{h(n)}{n^s},
\]
say. By our hypothesis $H(s)$ is bounded, and hence holomorphic, for $\si>1/2$. We modify the proof of Theorem 1 at several points. By Lemma 1 of \cite{Ka-Pe/2003} we have that for every $\epsilon>0$ there exists an integer $K=K(\epsilon)$ such that the coefficients $a^{-1}_G(n)$ of $1/G(s)$ satisfy
\[
a^{-1}_G(n) \ll n^\epsilon \hskip2cm (n,K)=1,
\]
and hence
\[
h(n)\ll n^{2\epsilon} \hskip2cm (n,K)=1.
\]
Therefore the set
\[
S = \{p: |h(p^m)|>p^{m/10} \ \text{for some $m\geq1$ or $p\leq 10^4$}\}
\]
is finite and we write
\begin{equation}
\label{1.13}
\begin{split}
H(s) &= \prod_p \frac{F_p(s)}{G_p(s)} = \prod_p H_p(s) \\
&=\prod_{p\not\in S}\left(1+\frac{h(p)}{p^s} + \frac{h(p^2)}{p^{2s}}\right) \prod_{p\in S}H_p(s)
\prod_{p\not\in S}\left(H_p(s)\left(1+\frac{h(p)}{p^s} + \frac{h(p^2)}{p^{2s}}\right)^{-1}\right)\\
&=Q_1(s)Q_2(s)Q_3(s),
\end{split}
\end{equation}
say. As in the prof of Theorem 1, $Q_2(s)$ and $1/Q_2(s)$ are holomorphic and bounded for $\si\geq1/2$. Moreover we have
\[
Q_3(s) = \prod_{p\not\in S}\left(1 + \frac{\sum_{m=3}^\infty \frac{h(p^m)}{p^{ms}}}{1+\frac{h(p)}{p^s} + \frac{h(p^2)}{p^{2s}}}\right) =  \prod_{p\not\in S}\left(1+ \sum_{m=3}^\infty \frac{k(p^m)}{p^{ms}}\right),
\]
say, and a computation shows that for $\si\geq 1/2$
\[
 \sum_{m=3}^\infty \frac{|k(p^m)|}{p^{m\si}} \leq \frac{1}{3} \ \  \text{for every $p\not\in S$ and} \ \  \sum_{p\not\in S} \sum_{m=3}^\infty \frac{|k(p^m)|}{p^{m\si}} \ll 1.
\]
Therefore, no factor of the product vanishes, and $Q_3(s)$ and $1/Q_3(s)$ are holomorphic and bounded for $\si\geq1/2$ as well.

\medskip
In order the treat $Q_1(s)$ we need the following elementary lemma.

\medskip
{\bf Lemma.} {\sl For every $a,b\in\CC$ there exists $\theta\in\CC$ with $|\theta|=1$ such that}
\[
|1+\theta a + \theta^2b| \geq 1 + \frac{1}{24}(|a|+|b|).
\]

\medskip
{\it Proof.} Suppose first that $|a|\leq |b|/2$. Then
\[
\max_{|\theta|=1} |1+\theta a + \theta^2b| \geq 1+|b|-|a| \geq 1+\frac12|b| \geq 1+\frac{1}{3}(|a|+|b|),
\]
and the result follows in this case. In the opposite case $|a|> |b|/2$ we apply the maximum 
\eject
modulus principle to the function $f(z) = 1 +az+bz^2$, thus obtaining
\[
\begin{split}
\max_{|\theta|=1} |1+\theta a + \theta^2b| &\geq \max_{|\theta|=1} |1+\frac{1}{4}\theta a +\frac{1}{16} \theta^2b| \\
&\geq 1+\frac{1}{4}|a| -\frac{1}{16}|b| \geq 1+ \frac{1}{24} (|a|+|b|),
\end{split}
\]
and the Lemma follows. Note that the constant $1/24$ is neither optimal nor important in what follows; moreover, in general it cannot be made arbitrarily close to 1. \fine

\medskip
From \eqref{1.13}, our hypothesis and the above information on $Q_2(s)$ and $Q_3(s)$ we deduce that there exists $M>0$ such that for $\si>1/2$
\[
|Q_1(s)| = \prod_{p\not\in S}\left|1+p^{-it}\frac{h(p)}{p^\si} +p^{-2it} \frac{h(p^2)}{p^{2\si}}\right| \leq M.
\]
By the Lemma, for every $\si$ and $p$ there exists $|\theta_{p,\si}|=1$ such that
\[
\left|1+\theta_{p,\si}\frac{h(p)}{p^\si} +\theta_{p,\si}^2\frac{h(p^2)}{p^{2\si}}\right| \geq 1 + \frac{1}{24}\left(\frac{|h(p)|}{p^\si} +\frac{|h(p^2)|}{p^{2\si}}\right).
\]
Assuming that $\si>1/2$ and $p\not\in S$, applying Kronecker's theorem as in the last part of the proof of Theorem 1 we find that
\[
\prod_{p\not\in S} \left(1 + \frac{1}{24}\left(\frac{|h(p)|}{p^\si} +\frac{|h(p^2)|}{p^{2\si}} \right)\right) \leq M.
\]
Then, letting $\si\to 1/2^+$, we deduce that the product
\[
\prod_{p\not\in S} \left(1 + \frac{1}{24}\left(\frac{|h(p)|}{p^{1/2}} +\frac{|h(p^2)|}{p} \right)\right)
\]
is convergent. Thus the series
\[
\sum_{p\not\in S} \left(\frac{|h(p)|}{p^{1/2}} +\frac{|h(p^2)|}{p} \right)
\]
is convergent as well and, in turn, the product
\[
\prod_{p\not\in S} \left(1 + \left(\frac{|h(p)|}{p^{1/2}} +\frac{|h(p^2)|}{p} \right)\right)
\]
converges. Hence $Q_1(s)$ and $Q_1(s)^{-1}$ are non-vanishing for $\si\geq 1/2$. 

\medskip
From \eqref{1.13} and the above properties of $Q_j(s)$, $j=1,2,3$, we immediately see that $H(s)$ is holomorphic and non-vanishing for $\si\geq 1/2$. Denoting by $\gamma_F(s)$ and $\gamma_G(s)$ the $\gamma$-factors of $F(s)$ and $G(s)$, thanks to the functional equation we deduce that
\[
\frac{\gamma_F(s)}{\gamma_G(s)}H(s)
\]
is a non-vanishing entire function of order $\leq 1$, and hence by Hadamard's theory we have
\begin{equation}
\label{1.14}
H(s) = \frac{\gamma_G(s)}{\gamma_F(s)} e^{as+b}
\end{equation}
with some $a,b\in\CC$. Now we can conclude by means of the almost periodicity argument that we used in our proof of the multiplicity one property of $\S$. For this we refer to Lemma 2.1 of \cite{Ka-Pe/1999b} and to Theorem 2.3.2 of \cite{Kac/2006}; in particular, \eqref{1.14} is exactly the last displayed formula of p.167 of \cite{Kac/2006}. This way we get that $H(s)\equiv1$, hence Theorem 3 is proved.

\medskip

\ifx\undefined\bysame{poly}.
\newcommand{\bysame}{\leavevmode\hbox to3em{\hrulefill}\ ,}
\fi

\medskip
\noindent
Jerzy Kaczorowski, Faculty of Mathematics and Computer Science, A.Mickiewicz University, 61-614 Pozna\'n, Poland and Institute of Mathematics of the Polish Academy of Sciences, 
00-956 Warsaw, Poland. e-mail: kjerzy@amu.edu.pl

\medskip
\noindent
Alberto Perelli, Dipartimento di Matematica, Universit\`a di Genova, via Dodecaneso 35, 16146 Genova, Italy. e-mail: perelli@dima.unige.it


\begin{thebibliography}{100} {\normalsize

\bibitem{Boh/1913} H.Bohr - {\sl \"Uber die Bedeutung der Potenzreihen unendlich vieler Variabeln in der Theorie der Dirichletschen Reihen $\sum\frac{a_n}{n^s}$} - Nachr. Akad. Wiss. G\"ottingen Math.-Phys. Kl. (1913), 441--488; see also {\sl Collected Mathematical Works}, vol.I, Dansk Matematisk Forening, 1952.

\bibitem{Bo-He/1995} E.Bombieri, D.A.Hejhal - {\sl On the distribution of zeros of linear combinations of Euler products} - Duke Math. J. {\bf 80} (1995), 821--862.

\bibitem{Br-He/2015} O.F.Brevig, W.Heap - {\sl Convergence abscissas for Dirichlet series with multiplicative
 coefficients} - to appear in Expo. Math. 2015; see arXiv:1506.04546v2.

\bibitem{Cha/1968} K.Chandrasekharan - {\sl Introduction to Analytic Number Theory} - Springer Verlag 1968.

\bibitem{Co-Gh/1993} J.B.Conrey, A.Ghosh - {\sl On the Selberg class of Dirichlet series: small degrees} - Duke Math. J. {\bf 72} (1993),  673--693.

\bibitem{Hec/1983} E.Hecke - {\sl Lectures on Dirichlet Series, Modular Functions and Quadratic Forms} - Vanderhoeck $\&$ Ruprecht 1983.

\bibitem{Kac/2006} J.Kaczorowski - {\sl Axiomatic theory of $L$-functions: the Selberg class} - In {\sl Analytic Number Theory}, C.I.M.E. Summer School, Cetraro (Italy) 2002, ed. by A.Perelli and C.Viola, 133--209, Springer L.N. 1891, 2006.

\bibitem{Ka-Pe/1999a} J.Kaczorowski, A.Perelli - {\sl On the structure of the Selberg class, I: $0\leq d \leq 1$} - Acta Math. {\bf 182} (1999), 207--241.

\bibitem{Ka-Pe/1999b} J.Kaczorowski, A.Perelli - {\sl The Selberg class: a survey} - In {\sl Number Theory in  Progress}, Proc. Conf. in Honor of A.Schinzel, ed. by K.Gy\"ory {\sl et al.},  953--992, de Gruyter 1999.

\bibitem{Ka-Pe/2003} J.Kaczorowski, A.Perelli - {\sl On the prime number theorem for the Selberg class} - Arch.
Math. {\bf 80} (2003), 255--263.

\bibitem{Ka-Pe/2005} J.Kaczorowski, A.Perelli - {\sl On the structure of the Selberg class, VI: non-linear twists} - Acta Arith. {\bf 116} (2005), 315--341.

\bibitem{Ka-Pe/book} J.Kaczorowski, A.Perelli - {\sl Introduction to the Selberg Class of $L$-Functions} - In preparation.

\bibitem{Lan/1915} E.Landau - {\sl \"Uber die Anzahl der Gitterpunkte in gewissen Bereichen (Zweite Abhandlung)} - Nachr. Ges. Wiss. G\"ottingen. Math.-Phys. Kl. (1915), 209--243.

\bibitem{Luo/1995} W.Luo - {\sl Zeros of Hecke $L$-functions associated with cusp forms} Acta Arith. {\bf 71} (1995), 139--157.

\bibitem{Ma-Qu/2010} B.Maurizi, H.Queff\'elec - {\sl Some remarks on the algebra of bounded Dirichlet series} -  J. Fourier Anal. Appl. {\bf 16} (2010), 676--692.

\bibitem{Ogg/1969} A.Ogg - {\sl Modular Forms and Dirichlet Series} - Benjamin 1969.

\bibitem{Per/2005} A.Perelli - {\sl A survey of the Selberg class of $L$-functions, part I} - Milan J. Math. {\bf 73} (2005), 19--52.

\bibitem{Per/2004} A.Perelli - {\sl A survey of the Selberg class of $L$-functions, part II} - Riv. Mat. Univ. Parma (7) {\bf 3*} (2004), 83--118.

\bibitem{Per/2010} A.Perelli - {\sl Non-linear twists of $L$-functions: a survey} - Milan J. Math. {\bf 78} (2010), 117--134.

\bibitem{Sel/1992} A.Selberg - {\sl Old and new conjectures and results about a class of Dirichlet series} - In {\sl Proc. Amalfi Conf. Analytic Number Theory}, ed by E.Bombieri {\it et al.}, 367--385, Universit\'a di Salerno 1992; {\sl Collected Papers} vol.II, 47--63, Springer 1991.

}
\end{thebibliography}
\end{document}